\documentclass[a4paper]{article}

\usepackage[T1]{fontenc}
\usepackage[utf8]{inputenc}
\usepackage{amsmath}
\usepackage{amsthm}
\usepackage{amssymb}
\usepackage{pgfplots}
\usepackage{relsize}
\usepackage{xcolor}
\usepackage{tikz}
\usetikzlibrary{plotmarks,mindmap,trees,shapes,arrows,shapes.multipart,patterns,fit,calc,positioning,decorations.pathreplacing,intersections,spy,through,backgrounds,decorations.markings}
\usepackage{vmargin}
\usepackage{MnSymbol}
\usepackage{enumitem}
\usepackage{multirow}
\usepackage[noend]{algpseudocode}
\usepackage{setspace}
\usepackage{arydshln}
\usepackage{caption}
\usepackage{adjustbox}
\usepackage{footmisc}
\usepackage[lofdepth,lotdepth]{subfig}
\usepackage[strict]{changepage}

\makeatletter
\newif\if@restonecol
\makeatother

\usepackage[algo2e,ruled,lined,titlenumbered]{algorithm2e}

\theoremstyle{remark}
\newtheorem{remark}{Remark}
\newtheorem{prop}{Property}

\pgfplotsset{compat=newest}

\newcommand{\s}{\ensuremath{^{(s)}}}

\newcommand{\sT}{\ensuremath{^{(s)^T}}}

\newcommand{\sj}{\ensuremath{^{(j)}}}
\newcommand{\sjT}{\ensuremath{^{(j)^T}}}
\newcommand{\invsj}{\ensuremath{^{(j)^{-1}}}}
\newcommand{\overj}{\ensuremath{^{(\overline{j})}}}
\newcommand{\overjT}{\ensuremath{^{(\overline{j})^T}}}
\newcommand{\invovj}{\ensuremath{^{(\overline{j})^{-1}}}}

\newcommand{\neij}{\ensuremath{^{\text{neigh}(j)}}}
\newcommand{\invneij}{\ensuremath{^{\text{neigh}(j)^{-1}}}}

\newcommand{\ddv}{\ensuremath{^{\diamondvert}}}
\newcommand{\ddm}{\ensuremath{^{\diamondminus}}}
\newcommand{\ddd}{\ensuremath{^{\diamondbackslash}}}

\newcommand{\ddmT}{\ensuremath{^{\diamondminus^T}}}
\newcommand{\dddT}{\ensuremath{^{\diamondbackslash^T}}}

\begin{document}

\title{A new impedance accounting for short and long range effects in mixed substructured formulations of nonlinear problems}

\author{Camille Negrello$^1$, Pierre Gosselet$^1$, Christian Rey$^2$ \\ \\ $^1$:LMT, ENS Paris Saclay/CNRS/Univ. Paris Saclay, \\61 avenue du Président Wilson, 94230 Cachan, France \\$^2$:Safran Tech, rue des Jeunes Bois, Chateaufort, \\CS 80112, 78772 Magny-les-Hameaux, France}

\date{}

\maketitle

\begin{abstract}
An efficient method for solving large nonlinear problems combines Newton solvers and Domain Decomposition Methods (DDM). In the DDM framework, the boundary conditions can be chosen to be primal, dual or mixed. The mixed approach presents the advantage to be eligible for the research of an optimal interface parameter (often called impedance) which can increase the convergence rate. The optimal value for this parameter is usally too expensive to be computed exactly in practice: an approximate version has to be sought, along with a compromise between efficiency and computational cost. In the context of parallel algorithms for solving nonlinear structural mechanical problems, we propose a new heuristic for the impedance which combines short and long range effects at a low computational cost.
\end{abstract}


\maketitle


\section{Introduction}

Dealing with nonlinear phenomena has become one of the predominant issues for mechanical engineers, in the objective of virtual testing. Whether they are geometrical or related to the material behavior, nonlinearities can be treated by a combination of Newton and linear solvers. Newton algorithms can be modified, secant, quasi-Newton \cite{crisfield1979faster,deuflhard1975relaxation,dennis1977quasi,zhang1982modified}, depending mostly on the complexity of tangent operators computation. If the meshed structure has a large number of degrees of freedom, linear solvers are chosen to be iterative and parallel, belonging to the class of Domain Decomposition Methods for instance \cite{mandel1993balancing,le1994domain,rixen1999simple,farhat2001feti,gosselet2006non}.

This article focuses on the nonlinear substructuring and condensation method, which has been investigated in previous studies \cite{cresta2007nonlinear,pebrel2008nonlinear,bordeu2009balancing,negrello2016substructured}. The substructured formulation involves a choice of interface transmission conditions, which can be either primal, dual or mixed, referring either to interface displacements, nodal interface reactions, or a linear combination of the two previous types -- i.e. Robin interface conditions. In this context, the mixed formulation has shown good efficiency \cite{cresta2008decomposition,pebrel2009etude,hinojosa2014domain,negrello2016substructured}, mostly due to a sound choice of the parameter introduced in the linear combination of interface conditions. Being homogeneous to a stiffness, and often refered to as an \textit{interface impedance}, this parameter can indeed be optimized, depending on the mechanical problem \cite{lions1990schwarz}. However, the computational cost of the optimal value involves in general storage and manipulation of global matrices, and is consequently not affordable in the framework of parallel computations. 

The interface impedance, in DDM methods for structural mechanics, should model, from the point of view of one  substructure, its interactions with the complement of the whole structure. In order to achieve good convergence rates without degrading computational speed, interface impedance can generally be approximated either by \textit{short scale} or \textit{long scale} formulations, depending on the predominant phenomena which must be accounted for. In the mechanical context, for instance, a common \textit{short scale} approximation can be built by assembling interface stiffness of the neighbors \cite{cresta2008decomposition, pebrel2009etude, hinojosa2014domain, negrello2016substructured}.
 
However, filtering long range interactions gives quite a coarse approximation of interface impedance, and does not give an accurate representation of the environment of each substructure. A good evaluation of the remainder of the structure should indeed couple these two strategies. Starting from this consideration, we propose here a new construction process of the interface impedance, based on a ``spring in series'' modeling of the structure, which couples the \textit{long} and \textit{short} range interactions with the structure. The heuristic we develop is strongly influenced by the availability of the various terms involved in our approximation. 

The first section of this paper introduces the reference (mechanical or thermal) problem and the notations used in the following. A succinct presentation of the nonlinear substructuring and condensation method then recalls the principles of its mixed formulation: how the interface nonlinear condensed problem is built from nonlinear local equilibriums, and the basics of the whole solving process, involving a global Newton algorithm combined with two internal solvers (parallel local Newton algorithms and a multi-scale linear preconditioned Krylov solver for tangent interface system). At Section~\ref{sec:new_heuristic}, the question of finding a relevant Robin parameter for mixed interface conditions is developed, mainly based on the observation that for each substructure, the optimal interface impedance is the nonlinear discretized \textit{Dirichlet-to-Neumann} operator of its complementary part. The new heuristic is then introduced, starting from the model of two springs in series, and a possible nonlinear multi-scale interpretation is given. The details of the two-scale approximation can be found at subsection \ref{ssec:two_scale}, its efficiency is evaluated at last section on several academic numerical examples.

\section{Reference problem, notations}

\subsection{Global nonlinear problem}

We consider here a nonlinear partial differential equation on a domain $\Omega$, representative of a structural mechanical or thermal problem, with Dirichlet conditions on the part $\partial \Omega_u \neq \varnothing$ of its boundary, and Neumann conditions on the complementary part $\partial \Omega_F$. After discretization with the Finite Element method, the problem to be solved reads:
\begin{equation}\label{eq:nl_problem}
f_{int}(u) + f_{ext} = 0
\end{equation}
Vector $f_{ext}$ takes into account boundary conditions (Dirichlet or Neumann) and dead loads, operator $f_{int}$ refers to the discretization of homogeneous partial differential equation.

\begin{remark} In linear elasticity, under the small perturbations hypothesis, one has:
\begin{equation*}
f_{int}(u) = -Ku
\end{equation*}
with $K$ the stiffness matrix of the structure.
\end{remark}

\subsection{Substructuring}

\begin{figure}[ht]
\begin{center}
\hspace{-0.5cm}\subfloat[Subdomains]{\includegraphics[width=4cm]{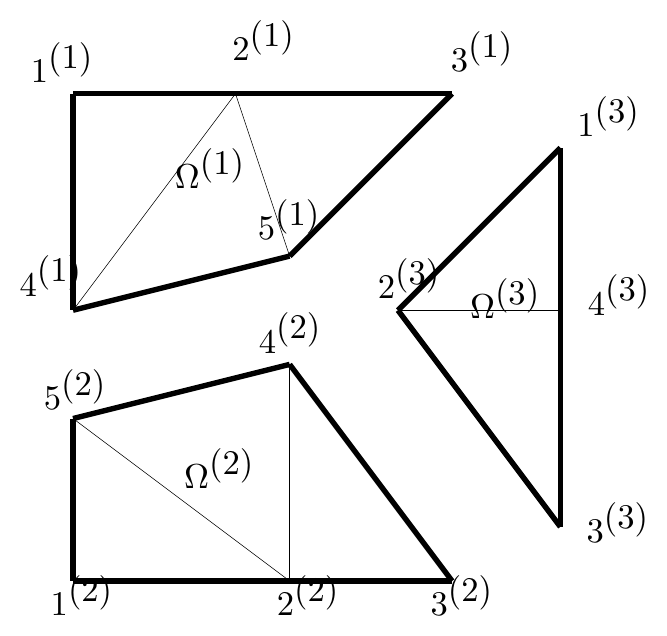}}
\subfloat[Local interface]{\includegraphics[width=4cm]{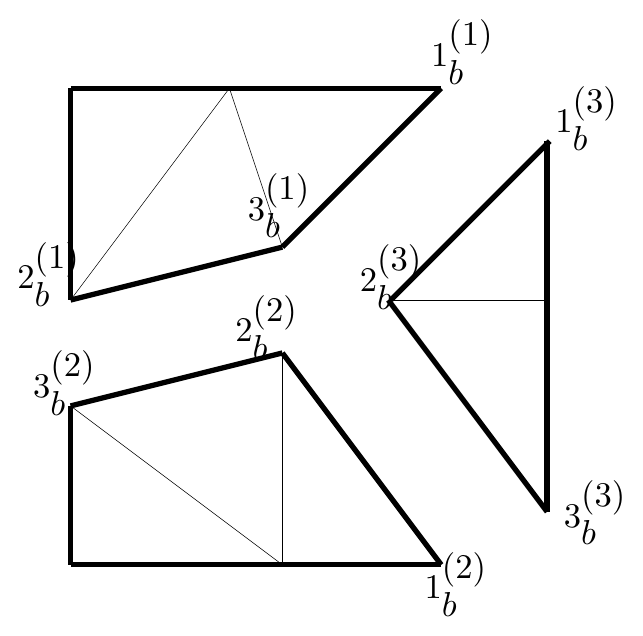}}
\subfloat[Interface nodes]{\includegraphics[width=3.5cm]{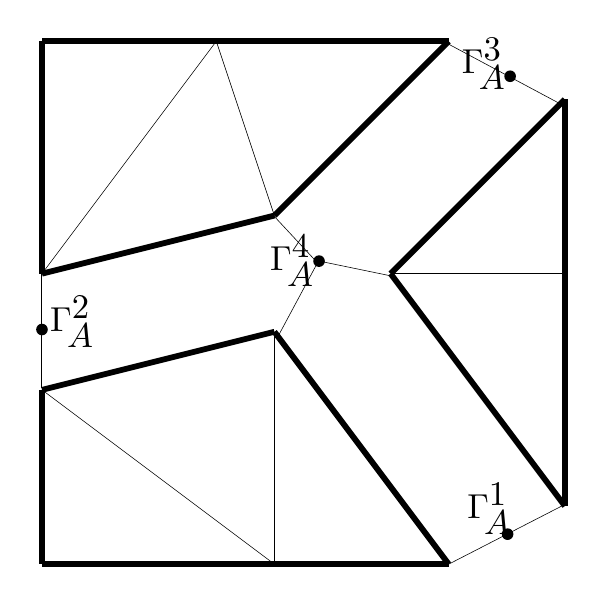}}
\subfloat[Interface connexions]{\includegraphics[width=3.5cm]{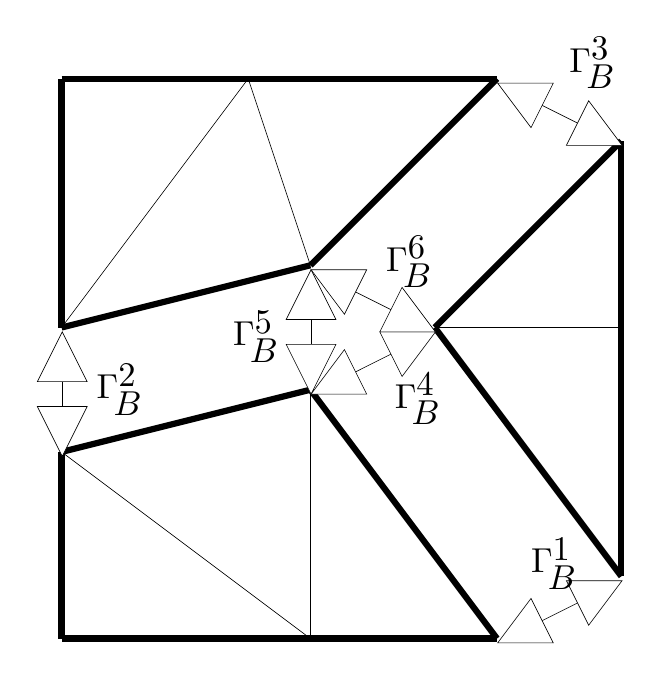}}

\hspace{-0.75cm}\subfloat{\includegraphics[width=14.5cm]{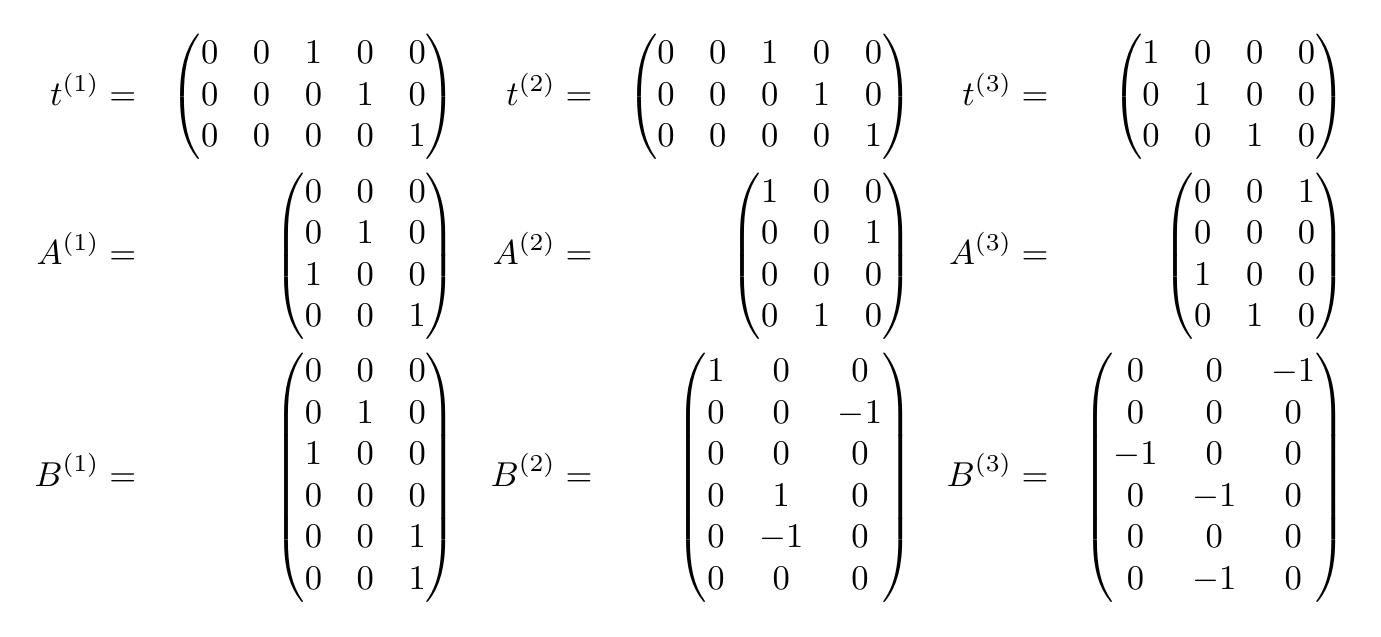}}
\end{center}
\caption{Local numberings, interface numberings, trace and assembly operators}
\label{fig:notations}
\end{figure}

Classical DDM notations will be used -- see figure \ref{fig:notations}: global domain $\Omega$ is partitioned into $N_s$ subdomains $\Omega\s$. For each subdomain, a trace operator $t\s$ restricts local quantities $x\s$ defined on $\Omega\s$ to boundary quantities $x\s_b$ defined on $\Gamma\s \equiv \partial \Omega\s \backslash \partial \Omega$:
\begin{equation*}
x\s_b = t\s u\s = x\s_{\vert \Gamma\s}
\end{equation*}
Quantities defined on internal nodes (belonging to $\Omega\s \backslash \Gamma\s$) are written with subscript $i$:  $x\s_i$.

Global primal (resp. dual) interface are noted $\Gamma_A$ (resp $\Gamma_B$).
Primal assembly operators $A\s$ are defined as canonical prolongation operators from $\Gamma\s$ to $\Gamma_A$: $A\s$ is a full-ranked boolean matrix of size $n_A \times n_b\s$ - where $n_A$ is the size of global primal interface $\Gamma_A$ and $n_b\s$ the number of boundary degrees of freedom belonging to subdomain $\Omega\s$. 

Diamond notations are used in the following: for a  domain $\Omega$ substructured in $N_s$ subdomains $\left( \Omega\s \right)$, concatenated local variables are superscripted $\ddv$, $\ddm$ or $\ddd$, depending on the alignment.
\begin{equation*}
\begin{aligned}
x\ddv = & \begin{pmatrix} 
x^{(1)} \\
\vdots \\
 x^{(N_{s})} \end{pmatrix},\qquad x\ddm = \begin{pmatrix}
{x^{(1)}} \, \ldots \, {x^{(N_{s})}} \end{pmatrix},\qquad
 M\ddd = \begin{pmatrix} 
M^{(1)} & 0   & 0 \\
0 & \ddots &  0  \\
0 & 0  & M^{(N_{s})} \\
\end{pmatrix} \\
\end{aligned}
\end{equation*}
Any matrix $B\s$ satisfying $\text{Range}(B\ddmT) = \text{Ker}(A\ddm)$ can be assigned to dual assembly operator -- see figure~\ref{fig:notations} for the most classical choice.

\section{Nonlinear substructuring and condensation: mixed formulation} 

This section recalls the principle of nonlinear substructuring and condensation, which is explained in details in \cite{negrello2016substructured}.

\subsection{Formulation of the condensed problem}

 Nonlinear problem \eqref{eq:nl_problem} is decomposed into $N_s$ nonlinear subproblems:
\begin{equation*}
f_{int}\ddv(u\ddv) + f_{ext}\ddv + t\dddT \lambda\ddv_b = 0\ddv
\end{equation*}
where $\lambda\s_b$ is the unknown local interface nodal reaction, introduced to represent interactions of the subdomain $\Omega\s$ with neighboring subdomains.

Transmission conditions hold:
\begin{equation*}
\left\lbrace \begin{aligned}
B\ddm u\ddv_b = 0 \\
A\ddm \lambda\ddv_b = 0
\end{aligned} \right.
\end{equation*}

The mixed formulation consists in introducing a new interface unknown:
\begin{equation*}
\mu\ddv_b = \lambda\ddv_b + Q\ddd_b u\ddv_b
\end{equation*}
where the matrix $Q\ddd_b$ is a parameter of the method. It has to be symmetric positive definite, and can be interpreted as a stiffness added to the interface, per subdomain: $Q\ddd_b$ is called \textit{interface impedance}. \medskip

Local equilibriums can then be reformulated as:
\begin{equation}\label{eq:local_eq}
f_{int}\ddv(u\ddv) + f_{ext}\ddv + t\dddT \left( \mu_b\ddv - Q\ddd_b u\ddv_b \right) = 0
\end{equation}
We assume the existence, at least locally, of a nonlinear mixed analogue $H\ddv_{nl}$ of the Schur complement (ie. a discrete Robin-to-Dirichlet operator):
\begin{equation}\label{eq:local_eq_H}
u\ddv_b = H\ddv_{nl}\left( \mu\ddv_b; Q\ddd_b, f\ddv_{ext} \right)
\end{equation}
\begin{prop} The tangent operator~$H\ddd_{t}$ to~$H\ddv_{nl}$ can be explicitly computed in function of the tangent stiffness~$K\ddd_t$:
\begin{equation*}
\begin{aligned}
 H\ddd_t = \frac{\partial H\ddv_{nl}}{\partial \mu\ddv_b} = t\ddd \left( K_t\ddd + t\dddT Q\ddd_b t\ddd \right)^{-1} t\dddT \\
\end{aligned}
\end{equation*}
Moreover, in the linear case, the Robin-to-Dirichlet operator written $H\ddv_{l}$ is  affine, with the constant term associated with external forces:
\begin{equation*}
\begin{aligned}
& H\ddv_l \left( \mu\ddv_b; Q\ddd_b, f\ddv_{ext} \right) = H\ddd_t \mu\ddv_b + b\ddv_m \\
\text{with }& b\ddv_m = t\ddd \left( K\ddd + t\dddT Q\ddd_b t\ddd \right)^{-1} f\ddv_{ext}
\end{aligned}
\end{equation*}
\end{prop}

\begin{remark} 
For the upcoming discussion, we will make use of the nonlinear primal Schur complement (Dirichlet-to-Neumann, noted $S_{nl}\s$) which is such that $\lambda_b\s=S_{nl}\s(u_b\s;f_{ext}\s)$. The tangent primal Schur complement can be computed from the tangent stiffness matrix:
\begin{equation*}
S_t\s=K\s_{t_{bb}}-K\s_{t_{bi}}{K\s_{t_{ii}}}^{-1}K\s_{t_{ib}}
\end{equation*}
and we have $H_t\ddd = (S\ddd_t + Q\ddd_b)^{-1}$. Note that the tangent dual Schur complement (Neumann-to-Dirichlet) can be written as ${S_t\s}^{\dagger} = t\s{K_t\s}^{\dagger}{t\s}^T$. In the linear case, the primal Schur complement is an affine operator with the constant term due to the external load:
\begin{equation*}
\begin{aligned}
& S\ddv_l \left( u\ddv_b; f\ddv_{ext} \right) = S\ddd_t u\ddv_b + b\ddv_p \\
\text{with }& b\ddv_p = f_{ext_b}\ddv - K\s_{t_{bi}}{K\s_{t_{ii}}}^{-1}f_{ext_i}\ddv = (S_t\ddd+Q_b\ddd)b_m\ddv
\end{aligned}
\end{equation*}
\end{remark}

Thanks to the  complementarity between balanced and continuous quantities, and to the symmetry positive definiteness of $Q\ddd_b$, any boundary displacement (defined independently on neighboring subdomains) can be split in a unique way into a continuous field belonging to $\text{Ker}\left( B\ddm \right)$ and a balanced field belonging to $\text{Ker}\left( A\ddm Q\ddd_b \right)$. Thus, the transmission conditions can be written in terms of $\mu\ddv_b$ and $u\ddv_b$, and gathered in a single equation:
\begin{equation*}
A\ddmT \left( A\ddm Q\ddd_b A\ddmT \right)^{-1} A\ddm \mu\ddv_b - u\ddv_b = 0
\end{equation*}

Finally, interface condensed problem reads:
\begin{equation}\label{eq:interf_pb}
R\ddv_b(\mu_b\ddv) \equiv A\ddmT \left( A\ddm Q\ddd_b A\ddmT \right)^{-1} A\ddm \mu\ddv_b - H\ddv_{nl} \left( \mu\ddv_b; Q\ddd_b, f\ddv_{ext} \right) = 0
\end{equation} 

\subsection{Solving strategy}

\subsubsection{Newton-Krylov algorithm}

Nonlinear substructuring and condensation results in applying a global Newton algorithm to interface problem \eqref{eq:interf_pb} instead of problem \eqref{eq:nl_problem}. Three steps are then involved in the solving process:
\begin{enumerate}[label=(\roman*)]
\item Local solutions of nonlinear equilibriums \eqref{eq:local_eq} are computed by applying local Newton algorithms.
\item The interface mixed residual is assembled.
\item The interface tangent problem is solved by a DDM solver.
\end{enumerate}
Newton global algorithm can be written, with previous notations:
\begin{equation*}
\left\lbrace \begin{aligned}
&\frac{\partial R\ddv_b}{\partial \mu_b\ddv} d\mu\ddv_b + R\ddv_b = 0 \\
&\mu\ddv_b  += d\mu\ddv_b
\end{aligned} \right.
\end{equation*}
Tangent problem then reads:
\begin{equation}\label{eq:tg_pb}
\left( A\ddmT \left( A\ddm Q\ddd_b A\ddmT \right)^{-1} A\ddm - H\ddd_t \right) d\mu\ddv_b = H\ddv_{nl} \left( \mu\ddv_b, Q\ddd_b, f\ddv_{ext} \right) - A\ddmT \left( A\ddm Q\ddd_b A\ddmT \right)^{-1} A\ddm \mu\ddv_b
\end{equation}

\subsubsection{Alternative formulation}\label{sec:altern_formul}
Tangent problem \eqref{eq:tg_pb} could be treated by a FETI-2LM solver \cite{roux2009feti}. An equivalent formulation of problem \eqref{eq:interf_pb} is also possible, where the boundary interface unknown $\mu\ddv_b$ is replaced by a couple of interface unknowns $\left( f_B, v_A \right)$, $f_B$ being a nodal reaction and $v_A$ an interface displacement. Couple $\left( f_B, v_A \right)$ is made unique by imposing the three following conditions:
\begin{itemize}[label=$\circ$]
\item $f_B$ is balanced
\item $v_A$ is continuous 
\item $\mu\ddv_b = B\ddmT f_B + Q\ddd_b A\ddmT v_A$
\end{itemize}
With this formulation, tangent problem is expressed by:
\begin{equation}\label{eq:pb_tg2}
\begin{aligned}
&\left(A\ddm S\ddd_t A\ddmT\right) dv_A  = A\ddm \left( Q\ddd_b + S\ddd_t \right) b\ddv_m \\
\text{with }&b\ddv_m = H\ddv_{nl} \left( \mu\ddv_b; Q\ddd_b, f\ddv_{ext} \right) - A\ddmT v_A
\end{aligned}
\end{equation}
Equation \eqref{eq:pb_tg2} has the exact form of a BDD \cite{mandel1993balancing} problem. It can thus conveniently be solved with usual preconditioner and coarse problem. The following quantities can then be deduced:
\begin{equation}\label{eq:recup_mu_u_lam}
\begin{aligned}
& d\mu_b\ddv = S_t\ddd A\ddmT dv_A - A\ddm \left( Q\ddd_b + S_t\ddd \right) b_m\ddv \\
& du\ddv = \left( K_t\ddd + t\dddT Q_b\ddd t\ddd \right)^{-1} t\dddT \left( A\ddm \left[ Q\ddd_b + S_t\ddd \right] b_m\ddv + d\mu_b\ddv \right) \\
& du_b\ddv = t\ddd du\ddv \\
& d\lambda_b\ddv = S_t\ddd du_b\ddv - A\ddm \left( Q\ddd_b + S_t\ddd \right) b_m\ddv = d\mu_b\ddv - Q\ddd_b du_b\ddv
\end{aligned}
\end{equation}

\subsubsection{Typical algorithm}

Algorithm~\ref{alg:robin-bdd} sums up the main steps of the method with the mixed nonlinear local problems and primal tangent solver. For simplicity reasons, only one load increment was considered.

As can be seen in this algorithm, several convergence thresholds are needed:
\begin{itemize}
\item Global convergence criterion $\varepsilon_{NG}$: since our approach is mixed, the criterion not only controls the quality of the subdomains balance (as in a standard Newton approach) but also the continuity of the interface displacement which is measured by an appropriate norm written $\|\cdot\|_B$. 
\item Local nonlinear thresholds $\varepsilon_{NL}\ddv$, which are associated with the Newton processes carried out independently on subdomains.
\item The global linear threshold of the domain decomposition (Krylov)  solver $\varepsilon_{K}$ (here BDD). 
\end{itemize}
The other parameters of the method are the initializations of the various iterative solvers and the choice of the impedance matrices $Q_b\ddd$. 

\begin{algorithm2e}[!ht]
\DontPrintSemicolon
\KwSty{Define:}\;
$r_{nl}^{m\diamondvert}(u\ddv,\mu_b\ddv) = f_{int}\ddv(u\ddv)-t\dddT Q_b\ddd t\ddd u\ddv + t\dddT \mu_b\ddv+ f\ddv_{ext}$\;
\BlankLine
\KwSty{Initialization:}\;
$(u_0\ddv,\lambda_{b_0}\ddv)$ such that $B\ddm t\ddd u_0\ddv=0$ and $A\ddm \lambda_{b_0}\ddv=0$\;
\KwSty{Set} $k=0$\;
\KwSty{Define} $\mu_{b_k}\ddv= \lambda_{b_k}\ddv + Q_b\ddd t\ddd u_k\ddv$\;
\While{$\| r_{nl}^{m\diamondvert}(u_k\ddv,\mu_{b_k}\ddv)\|+ \| B\ddm t\ddd u\ddv\|_{B}>\varepsilon_{NG}$ }{%
  \KwSty{Local nonlinear step}:  \;
   \KwSty{Set} $u_{k,0}\ddv=u_{k}\ddv$ and  $j=0$\; 
  \While{$\| r_{nl}^{m\diamondvert}(u_{k,j}\ddv,\mu_{b_k}\ddv)\|>\varepsilon\ddv_{NL}$}{
    $u_{k,j+1}\ddv=u_{k,j}\ddv-\left(K_{t_{k,j}}\ddd + t\dddT Q_b\ddd t\ddd \right)^{-1} r_{nl}^{m\diamondvert}( u_{k,j}\ddv, \mu_{b_k}\ddv)$ \;
    \KwSty{Set} $j=j+1$
  }
  \KwSty{Linear right-hand side}:\;
  $b_{m_k}\ddv = A\ddmT \left( A\ddm Q_b\ddd A\ddmT \right)^{-1} A\ddm \mu_{b_k}\ddv -t\ddd u_{k,j}\ddv $\;
  $b_{p_k}\ddv = (S_{t_{k,j}}\ddd + Q_b\ddd)b_{m_k}\ddv $\;
  \KwSty{Global linear step}:\;
  \KwSty{Set} $dv_A^{0}=0$ and  $i=0$\; 
  \While{$\|b\ddv_{p_k}-\left( A\ddm\; S_{t_{k,j}}\ddd A\ddmT \right) dv^i\|>\varepsilon_{K}$}{
    Make BDD iterations (index $i$)
  }
  \KwSty{Set} $u_{k+1}\ddv = u_k\ddv + du_k^{i\diamondvert}$ and $  \lambda_{b_{k+1}}\ddv = \lambda_{b_k}\ddv + d\lambda_{b_k}^{i\diamondvert}$ using \eqref{eq:recup_mu_u_lam}\;
  \KwSty{Set} $k=k+1$\;
}
\caption{Mixed nonlinear approach with BDD tangent solver}
\label{alg:robin-bdd}
\end{algorithm2e}

\section{New heuristic for the interface impedance}\label{sec:new_heuristic}

\subsection{Motivation}\label{ssec:motivation}

The parameter $Q\ddd_b$ is involved all along the solving, and a special care should be paid to its computation. 

\medskip

In order to frame the ideas, let us consider the Robin-Robin algorithm with stationary iteration, in the nonlinear case, with nonlinear impedance. Starting from the initial guess $\mu_b\ddv=0$, we have the iterations of Algorithm~\ref{alg:robin-robin}.

\SetNlSty{texttt}{(}{)}
\begin{algorithm2e}[ht]
	\DontPrintSemicolon
\nl Parallel solve: $S_{nl}\ddv(u_b\ddv;f\ddv_{ext})+Q_{nl}\ddv(u_b\ddv)= \mu_b\ddv$\;
\nl Parallel post-processing: $\lambda_b\ddv = S_{nl}\ddv(u_b\ddv;f\ddv_{ext}) = \mu_b\ddv-Q_{nl}\ddv(u_b\ddv)$\;
\nl Assembly: $\bar{u}_b\ddv = A\ddmT \tilde{A}\ddm u_b\ddv$, and $\bar{\lambda}_b\ddv = \left( I - \tilde{A}\ddmT A\ddm \right) \lambda_b\ddv$\;
\nl Parallel update of interface unknown: $\mu_b\ddv = Q_{nl}\ddv(  \bar{u}_b\ddv)  + \bar{\lambda}_b\ddv$\;
\caption{Robin-Robin stationary iteration}
\label{alg:robin-robin}
\end{algorithm2e}	

\medskip

Assembled quantities $\bar{u}_b\ddv$ and $\bar{\lambda}_b\ddv$ are defined such that the interface conditions can be written as:
\begin{equation}\label{eq:new_interf_cond}
u_b\ddv  = \bar{u}_b\ddv \text{ and }
\lambda_b\ddv  = \bar{\lambda}_b\ddv
\end{equation}
and we assume the nonlinear local operators $Q_{nl}\ddv$ to be such that the equivalence between \eqref{eq:new_interf_cond} and the following equation is ensured:
\begin{equation}\label{eq:intQnl}
\left(Q_{nl}\ddv(u_b\ddv) - Q_{nl}\ddv(\bar{u}_b\ddv)\right) + \left(\lambda_b\ddv-\bar{\lambda}_b\ddv\right)=0
\end{equation}

Considering a given subdomain $\Omega\sj$, and writing $\Omega\overj$ its complement, we can condense the whole problem on its interface; the boundary displacement $u_b\sj$ must then be the solution to:
\begin{equation*}
S_{nl}\sj(u_b\sj;f_{ext}\sj) + S_{nl}\overj(u_b\sj;f_{ext}\overj) =0
\end{equation*}
Comparing this equation with line~(1) of algorithm~\ref{alg:robin-robin}, one can see that, starting from a zero initial guess $\mu_b\sj=0$, the method converges in only one iteration with $Q_{nl}\sj = S_{nl}\overj$: the ideal impedance is the Dirichlet-to-Neumann operator of the complement. \medskip

In order to further discuss the problem, we now consider the linear case, and we recall that we have:  $S_{l}\sj(u_b\sj) = S_{t}\overj u_b\sj + b_p\overj$. In that case, the optimal impedance is thus an affine operator whose linear part (which we will write $Q_b\sj$ in agreement with the development of previous section) accounts for the stiffness of the complement domain, whereas the constant part accounts for the external load on the complement part. Note that another point of view is to use a strictly linear impedance together with a good (non-zero) initialization for $\mu_b\ddv$ which should account for the external load on the complement domain.

The construction of a good constant part for the impedance is usually realized, in the linear case, by the introduction of a well-chosen coarse problem; this is discussed in subsection~\ref{ssec:multi}. 
In the nonlinear case, building a coarse problem which would connect all subdomains during their inner Newton loop seems complex; more, it would break the independent computations. It looks simpler to rely on a good initialization in order to propagate the right-hand side: this can be done at low cost by easing accuracy constraints in the first inner Newton loop (adapting $\varepsilon_{NL}$), and then using the multiscale solver of the global linear step. Note that in \cite{klawonn2017new} a coarse problem is built for nonlinear versions of FETIDP and BDDC but, again, it mainly serves to find a good initialization before independent parallel nonlinear solves.

We now focus on the construction of $Q_b\ddd$, i.e. the linear part of the impedance. In the linear case,  one can show \cite{magoules2004optimal,gander2011optimal} that for a slab-wise  decomposition of the structure (or a tree-like decomposition, i.e. whose connectivity graph has no cycles), the setting $Q_b\sj = S_t\overj$ is optimal, in the sense that the convergence is reached in a maximum number of iterations equal to the number of subdomains (iterations are only needed to propagate the right-hand side). If the convenient coarse grid is added, convergence can be extremely fast. For an arbitrary decomposition, the optimality of such a setting can theoretically be lost, because of the unclear propagation rate of the right-hand side \cite{nier1998remarques,gander2011optimal}. However, the pertinence of this value still seems to be ongoing, especially being given the difficulty to define a more relevant setting for a matrix operator $Q_t\ddd$.

\medskip

Starting from these considerations, let us further analyze the terms of the following expression of the interface impedance for a given subdomain $\Omega\sj$:
\begin{equation}\label{eq:opti_qb_lin}
Q_t\sj = S_t\overj = K_{bb}\overj - K_{bi}\overj K_{ii}\invovj K_{ib}\overj
\end{equation}

The first term, $ K\overj_{bb}$, accounts for very local interactions. It is sparse, and exactly has the fill-in of matrix $K_{bb}\sj$. The second term, $K\overj_{bi} {K\overj_{ii}}^{-1} K\overj_{ib}$, accounts for long range interactions, it depends on the whole structure (geometry and material), and couples all degrees of freedom together via in-depth interactions. It is thus a full matrix; this property can be seen as the consequence of the pseudo-differentiability of the underlying Steklov-Poincaré operator of which the Schur complement is the discretization. It is important to note the minus sign: the short range part is very stiff and the global effects mitigate it. 

Obviously, formula~\eqref{eq:opti_qb_lin} is intractable in a distributed environment. However, different strategies have been investigated to compute  approximations at low cost -- see next subsection for a quick review.

\medskip
In the nonlinear context, the use of a linear impedance is of course non-optimal. Moreover, the best linear impedance probably resembles the Schur complement of the remainder of the subdomain in the final configuration, which is of course unknown a priori. Our aim is then to try to find a heuristic which gives an easy-to-compute approximation of the formula~\eqref{eq:opti_qb_lin} to be applied to the initial tangent stiffness.

\subsection{Quick review}\label{ssec:review}

The question of finding a good approximation of the Schur complement of a domain is at the core of mixed domain decomposition methods like optimized Schwarz methods \cite{gander2006osm} or the Latin method \cite{ladeveze2000micro}. Studies have proved that they needed to reproduce short-range effects (like local heterogeneity) but also structural effects (like the anisotropy induced by the slenderness of plate structures \cite{SAAVEDRA.2012.1}). When one wishes to choose an invariant scalar (or tensor in case of anisotropy) for each interface, it can be beneficial to use a coarse model for its estimation \cite{SAAVEDRA.2016.hal.1}. A possibility in order to better model short-range interaction between interface nodes is to use Ventcell conditions instead of simple Robin conditions \cite{Hoang2014}; this enables to recover the same sparsity for the impedance as for the stiffness of the subdomain. An extreme strategy is to use (scalar) Robin conditions on the Riesz' image of the normal flux leading to a fully populated impedance matrix \cite{DESMEURE.2011.3.1}. A more reasonable strategy is to use a strip approximation of the Schur complement \cite{magoules_algebraic_2006}, which can also be computed by adding elements to the subdomains \cite{Oumaziz2017}, in the spirit of restricted additive Schwarz methods \cite{cai1999restricted}.

From an algebraic point of view, short range approximation $ K\overj_{t_{bb}}$ (or even $ \operatorname{diag}(K\overj_{t_{bb}})$) is sometimes used for FETI's preconditioner \cite{Far94bis}, where it is called lumped approximation. 
Let $\text{neigh}(j)$ be the set of the neighbors of subdomain $j$, we have
\begin{equation}\label{eq:ktbb_neigh}
\text{lumped: }K_{t_{bb},l}\neij \equiv K_{t_{bb}}\overj = A\sjT \left( \sum_{s \in \text{neigh}(j)} A\s K\s_{t_{bb}} A\sT \right) A\sj
\end{equation}
or even:
\begin{equation}\label{eq:ktbb_neigh_diag}
\text{superlumped: }K_{t_{bb},sl}\neij \equiv \operatorname{diag}\left(K_{t_{bb}}\overj\right) = A\sjT \left( \sum_{s \in \text{neigh}(j)} A\s \operatorname{diag}\left(K\s_{t_{bb}}\right) A\sT \right) A\sj
\end{equation}
Being an assembly among a few subdomains of sparse block-diagonal matrices, this term is quite cheap to compute, and does not require any extra-computations, since local tangent stiffnesses are calculated anyway at each iteration of the solving process. The efficiency of the simple approximation \eqref{eq:ktbb_neigh} has been studied, in the context of nonlinear substructuring and condensation, in some research works \cite{cresta2007nonlinear, hinojosa2014domain, negrello2016substructured}, and has given good results when tested on rather homogeneous structures of standard shape. \medskip

In the domain decomposition framework for linear problems, long range interactions are taken into account thanks to the coarse grid problems \cite{Far94bis,mandel1993balancing,nouy03b}, which enables the method to comply with Saint-Venant's principle. These are closely related to projection techniques inspired by homogenization \cite{ibrahimbegovic2003strong, feyel2000fe, ladeveze2000micro, guidault2007two, guidault2008multiscale} in order to get low rank approximations. Let $U$ be an orthonormal basis of a well chosen subspace of displacements, the approximation can be written as:
\begin{equation}\label{eq:long_scale}
S\overj \simeq U (U^T S\overj U) U^T
\end{equation}
Saint-Venant's principle imposes $U$ to contain at least the rigid body motions of $\Omega\sj$, for computational efficiency it can be complemented by affine deformation modes or by displacements defined independently by interfaces.\medskip

However, if \textit{short range} approximations do not provide enough information to give a good representation of the faraway structure influence on a substructure $\Omega\sj$, neither do \textit{long range} approximation give a good estimation of the near-field response to a sollicitation. Besides, in the context of small displacements, a lack of precision on the close structure is more problematic than the filtering of long range interactions: predominant mechanical reactions usually come from nearby elements of the mesh.
 
The best strategy for $Q\sj_b$ would combine both \textit{short} and \textit{long} range formulations, however this version has not been much investigated yet. In particular it is not that easy to ensure the positivity of the impedance if the two approximations are computed independently. In \cite{gendre2011two} an expensive scale separation was introduced in the context of non intrusive global/local computations where $\Omega\overj$ was somehow available (which is not the case in our distributed framework). We propose here a new expression for parameter $Q\sj_b$, in the context of nonlinear substructuring and condensation with mixed interface conditions, which combines short and long scale formulations, at low computational cost.

\subsection{Spring in series model}\label{ssec:springs}

\begin{figure}[ht]
\begin{center}
\includegraphics{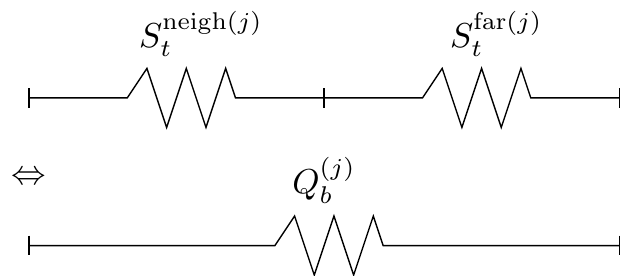}
\end{center}
\caption{Springs in series model}
\label{fig:springs_series}
\end{figure}

Our heuristic for the impedance relies on the simple observation that finding a two-scale approximation of the flexibility of $\Omega\overj$ may be more patent than for the stiffness. It is inspired by the simple model of two springs assembled in series: one spring models the stiffness of the neighboring subdomains whereas the second models the stiffness of the faraway subdomains (see figure \ref{fig:springs_series}). The resulting equivalent flexibility is the sum of the two flexibilities. 
In practice, in order to recover the structure of \eqref{eq:opti_qb_lin}, while remaining tractable, we propose the local flexibility ${S_t\neij}^{-1}$ to be the inverse of a sparse matrix, and the long-range flexibility ${S_t^{\text{far}(j)}}^{-1}$ to be low-rank. The latter condition is also motivated by \cite{bebendorf2003existence,amestoy2016complexity}, where it is shown that low-rank approximants of fully populated inverse operators, arising from FE discretization of elliptic problems, can be derived from the hierarchical-matrices theory. Typically we have:
\begin{equation}\label{eq:flex2terms}
\begin{aligned}
{Q_b\sj}^{-1} = {K_{t_{bb}}\neij}^{-1} + A\sjT V F V^T A\sj
\end{aligned}
\end{equation}
where $K_{t_{bb}}\neij$ can refer for instance to expressions \eqref{eq:ktbb_neigh} or \eqref{eq:ktbb_neigh_diag}, $F$ is a small-sized $m \times m$ square matrix, and $V$ an interface vectors basis of size $n_A \times m$. Writing $V\sj = A\sjT V$ the local contribution of basis $V$, expression \eqref{eq:flex2terms} can be inversed using the Sherman-Morrisson formula:
\begin{equation}
\begin{aligned}
Q_b\sj = K_{t_{bb}}\neij - \underset{W_b\sj}{\underbrace{ K_{t_{bb}}\neij V\sj}} \, \, \underset{ {M\sj}^{-1} }{\underbrace{ \left( F^{-1} + V\sjT K_{t_{bb}}\neij V\sj \right)^{-1}}} \, \, \underset{W_b\sjT}{\underbrace{ V\sjT K_{t_{bb}}\neij }}
\end{aligned}
\end{equation} 
This stiffness is a sparse matrix corrected by a low-rank term; then, when solving the (generalized) Robin problems, the Sherman-Morrisson formula can be used again:
\begin{equation}
\begin{aligned}
\text{let } \tilde{K}_t\sj & \equiv (K_t\sj+ t\sjT {K_{t_{bb}}^{\text{neigh}(j)} } t\sj)  
\text{ and } W\sj  \equiv t\sjT W_b\sj:\\
(K_t\sj+ t\sjT Q_b\sj t\sj)^{-1} & = \tilde{K}_t^{(j)^{-1}} + \tilde{K}_t^{(j)^{-1}} W\sj \left(M\sj - W\sjT \tilde{K}_t^{-1} W\sj \right)^{-1} W\sjT \tilde{K}_t^{(j)^{-1}}
\end{aligned}
\end{equation} 
The short-range term enables to regularize the problem without impairing the sparsity of the stiffness matrix.

\subsection{A multi-scale interpretation}\label{ssec:multi}

In the spirit of \cite{oumaziz2}, we can derive a multi-scale interpretation of the additive form \eqref{eq:flex2terms} adopted for the interface impedance. 

Starting from Algorithm~\ref{alg:robin-robin}, a macroscopic condition, inspired from the Latin method, can be imposed on the nodal reactions: the nodal reactions should satisfy a weak form of the interface balance, defined by a macroscopic basis $C_A$:
\begin{equation}\label{eq:macro_const}
C_A^T A\ddm \lambda_b\ddv = 0
\end{equation}
In the linear case, this condition can be enforced by the introduction of a Lagrange multiplier $\alpha$ (details can be found in \cite{oumaziz2}) in the interface condition \eqref{eq:intQnl}:
\begin{equation*}
\lambda_b\ddv - \bar{\lambda}_b\ddv + Q_b\ddd \left( u_b\ddv - \bar{u}_b\ddv \right) + Q_b\ddd A\ddmT C_A \alpha = 0
\end{equation*}
After algebraic calculations, writing local equilibriums with this new condition leads to:
\begin{equation*}
\left[ K\ddd + t\dddT Q_b\ddd \left( I\ddd - P_{C_A}\ddd \right) t\ddd \right] u\ddv = f_{ext}\ddv + t\dddT \left[ \bar{\lambda}_b\ddv + Q_b\ddd \left( I\ddd - P_{C_A}\ddd \right) \bar{u}_b\ddv \right]
\end{equation*}
where $P_{C_A}\ddd = A\ddmT C_A \left( C_A^T A\ddm Q_b\ddd A\ddmT C_A \right)^{-1} C_A^T A\ddm Q_b\ddd$ is a projector on the low-dimension subspace $\operatorname{Range}(A\ddmT C_A)$.

Not only the coarse space associated to the macroscopic constraint \eqref{eq:macro_const} results in the propagation of the right-hand side on the whole structure ($P_{C_A}\ddd$ is not sparse) but also in the modification of the impedance by the symmetric negative low rank term $-Q_b\ddd P_{C_A}\ddd$.

\medskip

In our nonlinear context, considering the basic setting $Q_b\sj = K_{t_{bb}}\neij$ \eqref{eq:ktbb_neigh} or \eqref{eq:ktbb_neigh_diag}, the modification $\bar{Q}_b\ddd = K_{t_{bb}}\neij - W_b\sj M\invsj W_b\sjT$ proposed in \eqref{eq:flex2terms} can be seen as the introduction of a multi-scale computation inside the mixed nonlinear substructuring and condensation method. As said earlier, the propagation of the right-hand side is ensured by a well-built initialization, which can be realized by adapting the inner Newton criterion $\epsilon_{NL}$ at each global iteration.

\subsection{Two-scale approximation of the flexibility}\label{ssec:two_scale}

\subsubsection{General idea}

From previous analysis, we try to derive an approximation of the (linear) optimal flexibility \eqref{eq:opti_qb_lin} which takes the additive form of \eqref{eq:flex2terms}. Being given a substructure $\Omega\sj$, we write $S_A\overj = \sum_{s \neq j} A\s S_t\s A\sT$ the assembly of local tangent Schur complements on the remainder $\Omega\overj$.

Using the quotient and the inverse formulas for the Schur complement, we have:
\begin{equation}\label{eq:flexSt}
{S_t\overj}^{-1} = \left( {S_A\overj}^{-1} \right)_{bb} = A\sjT {S_A\overj}^{-1} A\sj
\end{equation}

\begin{remark}
We here assume a substructuring ensuring the inversibility of $S_A\overj$ and $S_t\overj$, i.e. Dirichlet conditions are not concentrated on only one subdomain, and the complementary part of each subdomain is connected. In practice, this is almost always the case; if not, a simple subdivision can overcome the problem.
\end{remark}

\medskip

Classical preconditioners of BDD-algorithm can then be used as approximations of the inverse of $S_A\overj$. 
We hence introduce $\hat{G}_A\overj = \left[\, \ldots, \, \hat{A}\s_j R_b\s, \, \ldots \, \right]_{s \neq j}$ the concatenation of the scaled local traces of rigid body motions ($R_b\s$) of subdomains belonging to $\Omega\overj$, with $\hat{A}\s_j$ scaled assembly operators taking into account the absence of matter inside subdomain $\Omega\sj$. Considering the classical definition of scaled assembly operators $\tilde{A}\s$ \cite{klawonn2001feti}, modified operators $\hat{A}\s_j$ can be defined as:
\begin{equation*}
\begin{aligned}
 \hat{A}\s_j = & \left\lbrace \,\,\, \begin{aligned} \left( A\ddm \Delta\ddd A\ddmT - A\sj \Delta\sj A\sjT \right)^{-1} A\s \Delta\s \quad \text{ if } s \neq j \\
0 \quad \text{ if } s = j \end{aligned} \right. \\
& \hspace{0.5cm} \text{ with } \,\, \Delta\s \equiv \operatorname{diag} \left( K_{t_{bb}}\s \right)
\end{aligned}
\end{equation*}

Let $P_A\overj$ be the $S_A\overj$-orthogonal projector on $\operatorname{Ker}\left( \hat{G}_A\overjT S_A\overj \right)$:
\begin{align*}
 P_A\overj = I - \hat{G}_A\overj \left( \hat{G}_A\overjT S_A\overj \hat{G}_A\overj \right)^{-1} \hat{G}_A\overjT S_A\overj 
\end{align*}
we have:
\begin{equation}\label{eq:sep_inv_proj}
S_A\invovj =  P_A\overj S_A\invovj P_A\overjT + \left( I - P_A\overj \right) S_A\invovj \left( I - P_A\overj \right)^T 
\end{equation}
The BDD-theory states that in the first term, $S_A\invovj$ can be conveniently approximated by a scaled sum of local inverses\footnote{The GENEO theory \cite{SPILLANE:2013:FETI_GenEO_IJNME} states that, if needed, computable extra modes shall be inserted in $\hat{G}_A\overj$ in order to maintain the quality of the approximation.}.
After developing and factorizing, we have a first approximation of the flexibility:
\begin{equation}\label{eq:expr_inv_st}
\begin{aligned}
Q_{BDD}\invsj \equiv A\sjT \left( P_A\overj \sum_{s \neq j} \hat{A}\s_j {S_t\s}^{\dagger} \hat{A}\sT_j P_A\overjT +  \hat{G}_A\overj \left( \hat{G}_A\overjT S_A\overj \hat{G}_A\overj \right)^{-1} \hat{G}_A\overjT \right) A\sj \\
\end{aligned}
\end{equation}

\subsubsection{Long range interactions term}

The second term of expression \eqref{eq:expr_inv_st}, written $\hat{F}_{A,2}\sj$, is a matrix of low rank $m^{(j)}$, where  $m^{(j)}$ is the number of neighbors rigid body motions.  It could be used as is, however its computation involves the inversion of quantity $\hat{G}_A\overjT S_A\overj \hat{G}_A\overj$, an interface matrix of rank $m\overj$, where $m\overj$ is the number of local rigid body modes of the whole remainder $\Omega\overj$. In the context of large structures with a high number of subdomains, $m^{(\overline{j})}$ can increase drastically; saving the computation and factorization of such a matrix could then become quite interesting. Moreover, during the computation of the structure coarse problem, a close quantity is already assembled and factorized: the matrix $\tilde{G}_A^T S_A \tilde{G}_A$ -- with $S_A \equiv \sum_{s=1}^{N_s} A\s S\s A\sT$ and $\tilde{G}_A \equiv \left[ \ldots , \, \tilde{A}\s R_b\s, \, \ldots \right]$. Compared to $\hat{G}\overjT S_A\overj \hat{G}\overj$, the addition of the local term linked to $\Omega\sj$ in $\tilde{G}_A^T S_A \tilde{G}_A$ somewhat balances the classical scaling on its boundary (taking into account non-existant matter inside $\Omega\sj$), we thus propose:
\begin{equation*}
\hat{F}_{A,2}\sj \simeq A\sjT \hat{G}_A\overj \left( \tilde{G}_A^T S_A \tilde{G}_A \right)^{-1} \hat{G}_A\overjT A\sj \equiv \tilde{F}_{A,2}\sj
\end{equation*}

\subsubsection{Short range interactions term}

The first term of expression \eqref{eq:expr_inv_st}, written $\hat{F}_{A,1}\sj$, can also be simplified. First, for numerical efficiency, a diagonal lumping technique is used to approximate the local Schur complements (as explained in section \ref{ssec:review}). Then, in order to preserve sparsity, the projectors are removed. Assuming stiffness scaling is used we then directly recover the inverse of the superlumped stiffness of the neighbors:
\begin{equation}\label{eq:tild_ktbb}
\begin{aligned}
\hat{F}_{A,1}\sj & \simeq A\sjT \sum_{s \in \text{neigh}(j)} \hat{A}\s_j \operatorname{diag}\left( K_{t_{bb}}\s \right)^{-1} \hat{A}\sT_j A\sj \\ 
& = A\sjT \left( \sum_{s \in \text{neigh}(j)} A\s \operatorname{diag} \left( K_{t_{bb}}\s \right) A\sT \right)^{-1} A\sj = K_{t_{bb},\, sl}\invneij
\end{aligned}
\end{equation}

\subsubsection{Scaling issue}

A way to avoid building the modified scaled assembly operators $\hat{A}\s_j$ is to notice that for $s \neq j$, the following relation holds between modified and classical scaling operators $\tilde{A}\s$ \cite{klawonn2001feti}:
\begin{equation*}
\begin{aligned}
A\sjT \hat{A}_j\s & = \tilde{D}\sj A\sjT \tilde{A}\s \\
\text{with }\tilde{D}\sj \equiv A\sjT \left( A\ddm \Delta\ddd A\ddmT \right) & \left( A\ddm \Delta\ddd A\ddmT - A\sj \Delta\sj A\sjT \right)^{-1} A\sj
\end{aligned}
\end{equation*}
and we observe that the local diagonal matrix $\tilde{D}\sj$ can be 
extracted without cost from $\tilde{A}\sj$:
\begin{equation*}
\tilde{D}\sj = A\sjT \left( I - A\sj \tilde{A}\sjT \right)^{-1} A\sj
\end{equation*}
\begin{remark}
With evident notations, for a scaling based on the material stiffness, the diagonal coefficient of $\tilde{D}\sj$ associated with degree of freedom $x$ is equal to:
\begin{equation*}
\begin{aligned} 
\tilde{D}\sj_{xx} = \dfrac{ \sum_s K_{t_{xx}}\s }{ \sum_{s \neq j } K_{t_{xx}}\s } = \left( 1 - \dfrac{K_{t_{xx}}\sj }{ \sum_s K_{t_{xx}}\s} \right)^{-1} = \left( 1 - \tilde{A}\sj_x \right)^{-1}
\end{aligned}
\end{equation*}\qed
\end{remark}
\medskip

\noindent \textit{Final expression.} To conclude, we propose the following two-scale impedance:
\begin{equation}\label{eq:final_expr_qb}
\left( Q\sj_{b, \,2s} \right)^{-1} = K_{t_{bb},\,sl}\invneij + \tilde{D}\sj A\sjT \tilde{G}_A\overj \left( \tilde{G}_A^T S_A \tilde{G}_A \right)^{-1} \tilde{G}_A\overjT A\sj \tilde{D}\sj
\end{equation}

\subsection{Attempt to enrich the short-range approximation}\label{ssec:Qritz}

The short range part of the impedance, corresponding to the sparse approximation of $\hat{F}_{A,1}\sj$ by $K_{t_{bb},\,sl}\invneij$, seems very crude. In particular, we most probably underestimate the flexibility of the neighbors by using a diagonal operator.

We believe it is worth mentioning the tentative improvement which consisted in adding another low rank term:
\begin{equation*}
\hat{F}_{A,1}\sj \simeq K_{t_{bb},\,sl}\invneij + \tilde{D}\sj A\sjT V_k \Theta_k V_k^T A\sj \tilde{D}\sj
\end{equation*}
where $\Theta_k$ is a diagonal matrix and $V_k$ an orthonormal basis, approximations of the eigen-elements of $S_A\invovj$ associated with the higher part of the spectrum. They could be obtained at a moderate cost by post-processing the tangent BDD iterations in the spirit of \cite{gosselet2013total} (but considering the classical eigenvalues instead of the generalized ones). 

This low rank term could be concatenated with the one associated with rigid body motions $\tilde{F}_{A,2}\sj$, and thus did not modify the usability of the approximation.  
We observed that it led to a stiffness which was closer to our reference $S_t\invovj$ (measured with the Frobenius norm). But in practice when using it as the impedance in our numerical experiments, the reduction achieved in iterations numbers was not worth the additional cost of the enrichment term -- this is why we do not present it in detail. This ``improvement'' may be more useful on other classes of nonlinear problems for which it would be important not to overestimate the stiffness of the remainder of the structure.

\section{Results}

\subsection{Two test cases}\label{sec:two_test_cases}

\begin{figure}[!ht]
\includegraphics[width=\textwidth]{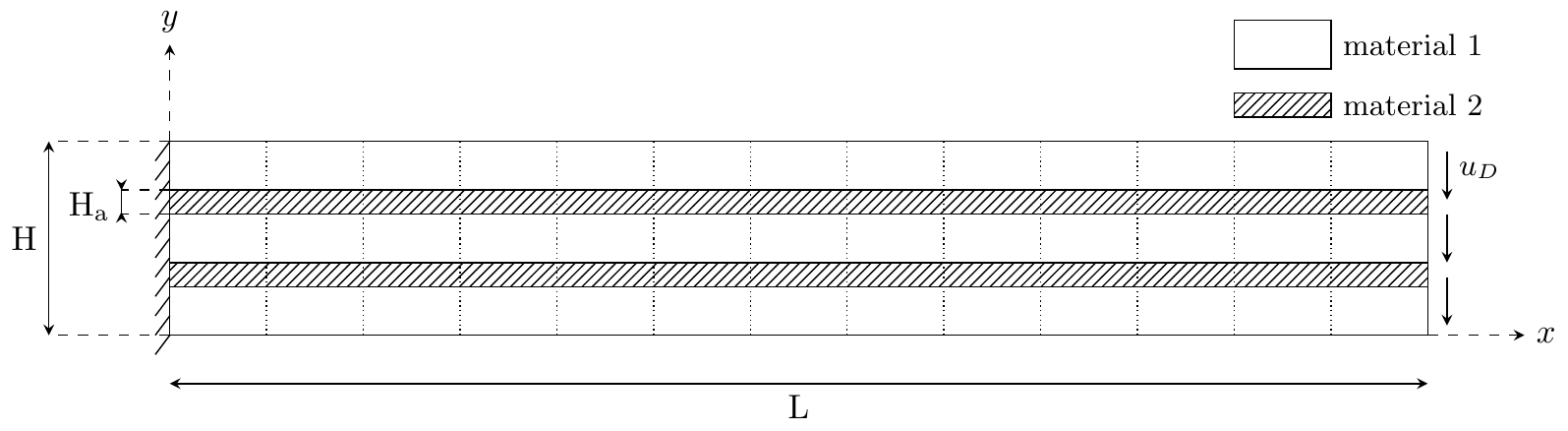}
\caption{Bi-material beam: partition and loading}
\label{fig:test_case_bimat}
\end{figure}

The efficiency of the expression \eqref{eq:final_expr_qb} is evaluated on two numerical test cases. First test case is a bi-material beam with bending load, represented on figure~\ref{fig:test_case_bimat}. Material and geometrical parameters are given in table \ref{tab:params_tests_case}: one of the two materials is chosen to be elastoplastic with linear hardening, the other one is chosen to remain elastic. Load is applied with imposed displacement on the edge defined by $x=L$.

\begin{figure}[!ht]
\begin{center}
\includegraphics[width=15cm]{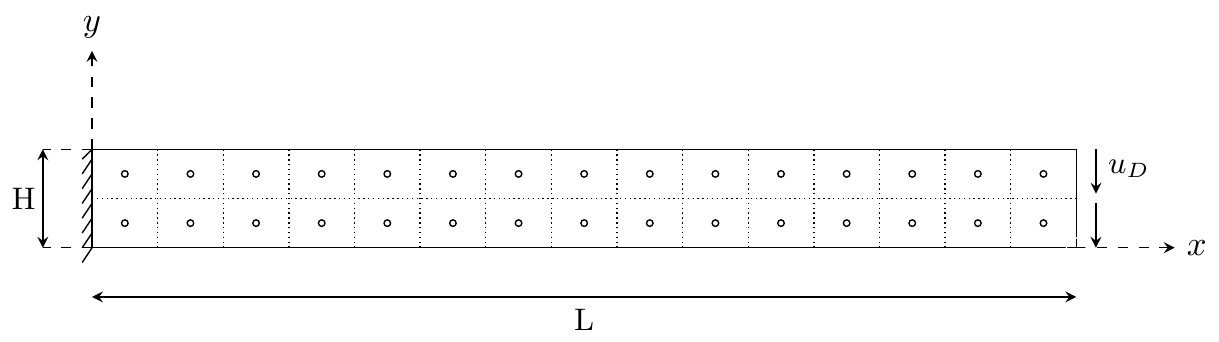}
\end{center}
\caption{Multiperforated beam: partition and loading}
\label{fig:test_case_multiperf}
\end{figure}

Second test case is a homogeneous multiperforated beam with bending load, represented on figure~\ref{fig:test_case_multiperf}. Material and geometrical parameters are given in table \ref{tab:params_tests_case}: material is chosen to be  elastoplastic with linear hardening. Load is applied with imposed displacement $u_D$ on the edge defined by $x=L$.

\begin{table}[ht]
\begin{center}
\begin{minipage}{0.56\linewidth}
\begin{tabular}{|c:c|c|}
\hline \hline
\multicolumn{3}{|c|}{Bi-material beam} \\
\hline \hline
\multicolumn{3}{|c|}{ Material parameters} \\
\hline \hline
 & Material 1 & Material 2 \\
Young & $E_1 = 420e2$ & $E_2 = 210e6$ \\
Poisson coefficient & $\nu_1 = 0.3$ & $\nu_2 = 0.3$ \\
Elastic limit & & $\sigma_{0_2} = 420e3$ \\
Hardening coefficient & & $h_2 = 1e3$   \\
\hline \hline
\multicolumn{3}{|c|}{Geometrical parameters} \\
\hline \hline
Total length & \multicolumn{2}{:c|}{L $ = 13$} \\
Total height & \multicolumn{2}{:c|}{H $ = 2$} \\
Height of an armature & \multicolumn{2}{:c|}{H$_\text{a} = 0.25$} \\
\hline \hline
\end{tabular}
\end{minipage}
\begin{minipage}{0.43\linewidth}
\begin{tabular}{|c:c|}
\hline \hline
\multicolumn{2}{|c|}{Multiperforated beam} \\
\hline \hline
\multicolumn{2}{|c|}{ Material parameters} \\
\hline \hline
 & \\
Young & $E = 210e6$ \\
Poisson coefficient &  $\nu = 0.3$ \\
Elastic limit & $\sigma_0 = 420e3$ \\
Hardening coefficient & $h = 1e6$ \\
\hline \hline
\multicolumn{2}{|c|}{Geometrical parameters} \\
\hline \hline
Length & L $ = 10$ \\
Height & H $ = 1$ \\
Hole radius & r $ = 2/30$ \\
\hline \hline
\end{tabular}
\end{minipage}
\caption{Material and geometrical parameters}
\label{tab:params_tests_case}
\end{center}
\end{table}

\subsection{Elastic analysis}

The ultimate goal of this paper is to assess the performance of the new impedance \eqref{eq:final_expr_qb} in the nonlinear multi-scale distributed context. Before we reach that point, a preliminary mono-scale elastic study is performed in order to verify that the heuristic developed in previous sections is actually able to capture both short and long range interactions within the structure.

Sollicitations are here keeped low enough to remain in the elastic domain of every materials: bi-material beam and multiperforated beam are both submitted to a bending load of intensity $u_D = 1.5 \,\, 10^{-3}$. More, decomposition is for now only performed along $x$-axis (multiple points will be involved in next section, where the nonlinear multi-scale context is considered). One of the interest of the elastic linear case with slab-wise decomposition relies on the ability to express the optimal interface impedance: $Q_b\sj = S_t\overj$ (see \ref{ssec:motivation}). Even if the computational cost of this parameter would be, in a real situation, absolutely not affordable in the context of parallel resolutions, it was calculated here for the purpose of our analysis. A comparison with an optimal reference can thus be made for the two following expressions:
\begin{itemize}[label=$\circ$]
\item a classical choice $K_{bb,l}\neij$: see \eqref{eq:ktbb_neigh}
\item the new expression $Q_{b,\,2s}\sj$: see \eqref{eq:final_expr_qb}
\end{itemize}
Being given the alternative formulation we chose for the mixed nonlinear substructuring and condensation method (see section \ref{sec:altern_formul}), an elastic resolution would be strictly equivalent to a primal BDD resolution. Therefore, no comparison of different interface impedances is possible with Algorithm~\ref{alg:robin-bdd}. A mono-scale FETI-2LM solver \cite{roux2009feti} was hence implemented, corresponding to the first formulation of the mixed interface problem with the $\mu_b\ddv$ unknown \eqref{eq:tg_pb}. This algorithm enables to solve linear problems with Robin interface transmission conditions. 

Note that an optimal coarse problem could be added in order to recover an efficient multi-scale solver \cite{dubois2012optimized,haferssas2015robust,loisel2015optimized}. However, this augmentation strategy would make it impossible to discern the efficiency of the long range interactions term of our two-scale impedance. Again, our aim is not to compete with augmented Krylov solvers for linear problems but to find an alternative way, compatible with nonlinear problems, to introduce long-range effects. The mono-scale formulation is thus preserved in order to evaluate the  ability of \eqref{eq:final_expr_qb} to introduce in local equilibriums information related to the interactions with the far structure, in a linear context where the optimal parameter is known.

\begin{table}[!ht]
\begin{center}
\subfloat[Bi-material beam]{\includegraphics[width=10cm]{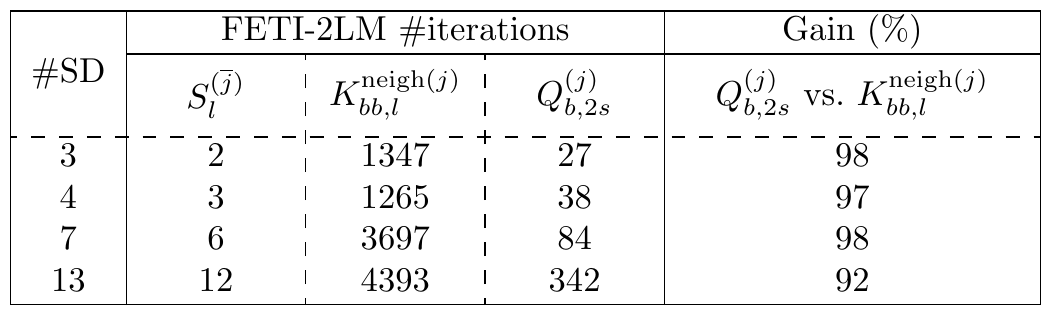}\label{subtab:results_lin_bimat}} \qquad
\subfloat[Multiperforated beam]{\includegraphics[width=10cm]{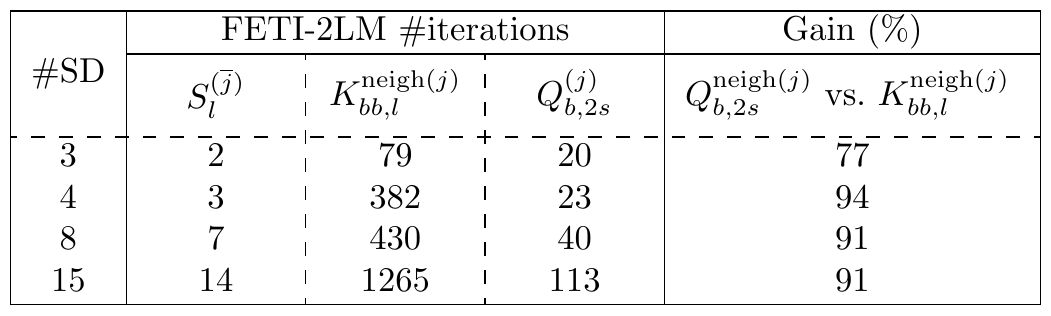}\label{subtab:results_lin_multiperf}}
\end{center}
\caption{Comparison of the three interface impedances: linear behavior}
\label{tab:results_lin}
\end{table}

Results are given on table \ref{tab:results_lin} for the two previously introduced test cases.

As expected, for both test cases, the optimal interface impedance $S_t^{(\overline{j})}$ rounds off the resolution after a number of iterations equal to the number of subdomains minus one. Being given the repartition of the subdomains (no multiple points) and the absence of a coarse problem, this is the best convergence rate that can be achieved: mixed transmission conditions with interface impedance $S_t^{(\overline{j})}$ is optimal.

The classical choice $K_{bb,l}^{\text{neigh}(j)}$ does not involve any information on the long range interactions of a subdomain with the faraway structure inside local equilibriums: the number of iterations drastically increases along with the number of substructures. 

The new expression $Q_{b,2s}^{(j)}$ introduced in this paper highly reduces the numbers of FETI-2LM iterations, compared to classical choice $K_{bb,l}^{\text{neigh}(j)}$: gains are between 77 and 98\%. This should mostly be due to the additive form of expression \eqref{eq:final_expr_qb}, with the introduction of a long range interactions-term in the flexibility -- obviously, the absence of coarse problem in the resolution reinforces the benefits of this term. The forthcoming nonlinear study, based on algorithm \ref{alg:robin-bdd}, will replace this expression in a context of multiscale computation. 

Performance of expression $Q_{b,2s}^{(j)}$ is evidently not as good as that of optimal expression $S_t^{(\overline{j})}$, but the increase in iterations numbers is only of about ten times the optimal iterations number (while it reaches about hundreds times the optimal iterations number for $K_{bb,l}^{\text{neigh}(j)}$). We also recall that interface impedance $S_t^{(\overline{j})}$ can not be computed in parallel resolutions: expression $Q_{b,2s}^{(j)}$, at the contrary, is fully and easily tractable.

The expression introduced here to evaluate the interface impedance thus seems, at least in the linear case, to achieve great performance at very low cost. 

\begin{remark} As said earlier, the second effect of multiscale approaches (beside modifying the Robin condition), lies in the instantaneous propagation of the right-hand side. In our approach, the absence of a coarse problem is somehow compensated by the presence of tangent interface systems (solved with state of the start multi-scale BDD method). As an example, we initialized our linear FETI-2LM solver with the fields resulting from one BDD iteration. For the multiperforated beam split in 15 subdomains, the number of FETI-2LM iterations goes from 113 to 89, which is significant (for less subdomains, the coarse problem is too small to bring any valuable piece of information). In the spirit of \cite{negrello2016substructured}, a tuned setting of the solvers' thresholds (synchronized with the evolution of the global residual, i.e. the precision of the global solution) could perform a good compromise between a global spread of the information and independent computations.

\end{remark}

\subsection{Plastic analysis}\label{sec:it_numbers}

The evaluation of the performance of expression \eqref{eq:final_expr_qb}  is continued with a plastic evolution study. The two test cases are submitted to bending loads, applied incrementally. Bi-material beam loading is decomposed as follows:
\begin{align}
u_D & = \left[ 0.05, \, \, 0.1, \, \, 0.15, \, \, 0.2, \, \, 0.25, \, \, 0.3, \, \, 0.35, \, \, 0.375, \, \, 0.4, \, \, 0.425, \, \, 0.45 \right] u_{max} \label{eq:load_incre_bimat}\\
u_{max} & = 7.1 \nonumber
\end{align}
For multiperforated beam loading, the incremental decomposition is set to:
\begin{align}
u_D & = \left[ 0.4, \, \, 0.6, \, \, 0.8, \, \, 1, \, \, 1.15, \, \, 1.3, \, \, 1.45, \, \, 1.5 \right] u_{max} \label{eq:load_incre_multiperf} \\
u_{max} & = 0.275 \nonumber
\end{align}
\begin{remark} 
For the sake of clarity, every over load increment \eqref{eq:load_incre_bimat} and \eqref{eq:load_incre_multiperf} is represented in the forthcoming results tables. 
\end{remark}

The substructuring of the bi-material beam involves 13 subdomains along $x$-axis, while multiperforated beam is decomposed into 30 subdomains with multiple points (see figures~\ref{fig:test_case_bimat} and~\ref{fig:test_case_multiperf}).

Numbers of Krylov iterations, cumulated over global Newton loops and load increments, are stored for the three interface impedances $S_t^{(\overline{j})}$, $K_{bb,l}^{\text{neigh}(j)}$ and $Q_{b,2s}^{(j)}$ and the two test cases in tables \ref{tab:iter_bimat} and \ref{tab:iter_multiperf}. Indeed, performance of the solver is in particular linked to the number of processor communications, which are directly proportional to the number of Krylov iterations. 

The computation of local tangent operators, at each global iteration, is also a costly operation. The numbers of global Newton iterations, cumulated over load increments, are thus also stored for each expression of the interface impedance and the two test cases. Note that in these cases, the number of Krylov iterations is almost constant per linear system, the cumulated numbers of Krylov iterations are thus nearly proportional to the numbers of global Newton iterations; the latter are therefore only stored for the last load increment. 

A fourth approach has been added to the study, written NKS in both tables, and corresponding to the ``classical'' resolution process used in nonlinear structural mechanical problems: a global Newton algorithm, combined with a linear DD solver for the tangent systems. The main difference between the nonlinear substructuring and condensation method and this classical technique resides in the nonlinear/linear algorithms used for local resolutions. The resulting comparisons with approaches $S_t^{(\overline{j})}$, $K_{bb,l}^{\text{neigh}(j)}$ and $Q_{b,2s}^{(j)}$ hence represent the gains that can be achieved with the mixed nonlinear substructuring and condensation method, in the more general framework of nonlinear solvers. \medskip

\begin{table}[!ht]
\begin{center}
\begin{minipage}{\linewidth}
\begin{center}
\begin{tabular}{|c|cccccc|c|}
\hline
\multicolumn{7}{|c|}{Krylov} & Global Newton \\
\hline
load inc. & 0.05 & 0.15 & 0.25 & 0.35 & 0.4 & 0.45 & 0.45\\
\hline
$S_t^{(\overline{j})}$ & 37 & 229 & x &  &  &  & x \\
$K_{bb,l}^{\text{neigh}(j)}$ & 36 & 259 & 598 & 978 & 1357 & 1772 & 47 \\
$Q_{b,2s}^{(j)}$ & 37 & 266 & 580 & 891 & 1204 & 1514 & 39 \\
\hline
\hline
NKS & 74 & 296 & 633 & 970 & 1344 & 1795 & 48 \\
\hline
\multicolumn{8}{c}{ } \\
\hline
\multicolumn{8}{|c|}{Gains (\%) } \\
\hline
$Q_{b,2s}^{(j)}$ vs. $K_{bb,l}^{\text{neigh}(j)}$ & -3 & -3 & 3 & 9 & 11 & 15 & 17 \\
\hline
\hline
$Q_{b,2s}^{(j)}$ vs. NKS & 50 & 10 & 8 & 8 & 10 & 16 & 19 \\
\hline
\end{tabular}
\end{center}
\end{minipage}
\end{center}
\caption{Bi-material beam: Krylov cumulated iterations over load increments, global Newton cumulated iterations}
\label{tab:iter_bimat}
\end{table}

\begin{table}[!ht]
\begin{center}
\begin{minipage}{\linewidth}
\begin{center}
\begin{tabular}{|c|cccc|c|}
\hline
\multicolumn{5}{|c|}{Krylov} & Global Newton\\
\hline
load inc. & 0.6 & 1 & 1.3 & 1.5 & 1.5\\
\hline
$S_t^{(\overline{j})}$ & 48 & 172 & 322 & 481 & 30 \\
$K_{bb,l}^{\text{neigh}(j)}$ & 61 & 212 & 373 & 548 & 35 \\
$Q_{b,2s}^{(j)}$ & 48 & 170 & 300 & 438 & 28 \\
\hline
\hline
NKS & 73 & 222 & 385 & 561 & 36 \\
\hline
\multicolumn{6}{c}{ } \\
\hline
\multicolumn{6}{|c|}{Gains (\%)} \\
\hline
$Q_{b,2s}^{(j)}$ vs. $K_{bb,l}^{\text{neigh}(j)}$ & 21 & 20 & 20 & 20 & 20 \\
$Q_{b,2s}^{(j)}$ vs. $S_t^{(\overline{j})}$ & 0 & 1 & 7 & 9 & 7 \\
\hline
\hline
$Q_{b,2s}^{(j)}$ vs. NKS & 34 & 23 & 22 & 22 & 22 \\
\hline
\end{tabular}
\end{center}
\end{minipage}
\end{center}
\caption{Multiperforated beam: Krylov cumulated iterations over load increments, global Newton cumulated iterations}
\label{tab:iter_multiperf}
\end{table}

A first preliminary observation compares results for interface impedance $S_t^{(\overline{j})}$ in the linear and the nonlinear case: the primitive guess we made about $S_t^{(\overline{j})}$ being the best possible approximation we could analytically define of the interface impedance value was mistaken in the nonlinear formulation. For bi-material beam for instance, the resolution ended up with a divergence in the local Newton solvers, caused by fake high levels of plasticity inside subdomains, artifacts of the resolution -- this may be due to an excessively soft interface impedance, which lets the material deform more than necessary. 

Secondly, although our first guess was apparently misguided, the additive expression we derived from it seems to behave very satisfyingly: best performance is now achieved -- in the nonlinear process -- with the new expression of interface impedance $Q_{b,2s}^{(j)}$. Gains in terms of Krylov cumulated iterations, compared to classical interface impedance $K_{bb,l}^{\text{neigh}(j)}$, vary from 15\% to 20\% at the end of the resolution: a benefit which should represent a non negligible decrease in CPU time for large structure problems (where each communication operation can be highly time-consuming). Compared to the interface impedance $S_t^{(\overline{j})}$ -- which is not computationally affordable in practice, -- only the multiperforated beam can be effectively studied (convergence was not reached for bi-material beam): gains, for approach $Q_{b,2s}^{(j)}$, reach up 9\% at the end of the resolution -- in terms of cumulated numbers of Krylov iterations. 

\begin{remark}
Bi-material beam was meshed with 25~789 degrees of freedom, and its substructuring into 13 subdomains involved 984 interface degrees of freedom. Multiperforated beam was meshed with 30~515 degrees of freedom, and its substructuring into 30 subdomains involved 1641 interface degrees of freedom. Despite the relative smallness of these test cases, we expect them to be representative of computations on larger structures. Unfortunately our Octave-based code did not allow meaningful time measurements and large scale computations. Moreover, limiting communication as we try to do would be even more appreciable on computations involving many processors. The number of Krylov iterations seems to be the fairest and most reliable performance measurement.
\end{remark}

Comparison with classic method shows similar results for both test cases: at the end of the resolution, gains vary from 16 to 22\% for Krylov cumulated iterations, and from 19 to 22\% for global Newton cumulated iterations. This gain corresponds to the overall performance of the nonlinear substructuring and condensation method that can be achieved with mixed approach, compared to classical procedures.

\begin{remark}
The rather limited performance of mixed nonlinear substructuring and condensation method with classical interface impedance $K_{bb,l}^{\text{neigh}(j)}$, compared to the classical resolution method, can be noticed in the above two examples. This lack of efficiency can probably be imputed to the difficulty of giving full account of long range phenomena with a short-scale interface impedance, whereas they prevail in the case of local heterogeneity (bi-material beam) and slenderness of plate structures (multiperforated beam).
\end{remark}

\subsection{Coupling with SRKS-method}

An augmentation strategy of Krylov subspaces, at each global nonlinear iteration, is possible by extracting Ritz vectors and values at the end of each Krylov solving and re-using them to construct an augmentation basis for the following Krylov iterations. The so-called TRKS method \cite{gosselet2013total} reuses all of the produced Ritz vectors, while SRKS method \cite{gosselet2013total} consists in selecting the Ritz values which are good enough approximations of tangent operator eigenvalues, and the corresponding Ritz vectors. SRKS method was implemented and its coupling with nonlinear substructuring and condensation method was studied for both test cases defined at section \ref{sec:two_test_cases}. 

Results are given in tables \ref{tab:iter_bimat_SRKS} and \ref{tab:iter_multiperf_SRKS}.

\begin{table}[!ht]
\begin{center}
\begin{tabular}{|c|cccccc:c|}
\hline
\multicolumn{8}{|c|}{Krylov}\\
\hline
\multicolumn{7}{|c:}{ with SRKS } & wo SRKS \\
\hline
load inc. & 0.05 & 0.15 & 0.25 & 0.35 & 0.4 & 0.45 & 0.45 \\
\hline
$S_t^{(\overline{j})}$ & 37 & 97 & x &  &  &  & x \\
$K_{bb,l}^{\text{neigh}(j)}$ & 36 & 113 & 218 & 339 & 452 & 578 & 1772 \\ 
$Q_{b,2s}^{(j)}$ & 37 & 98 & 197 & 304 & 398 & 492 & 1514 \\
\hline
\hline
NKS & 53 & 130 & 245 & 370 & 492 & 648 & 1795 \\
\hline
\multicolumn{8}{c}{ } \\
\hline
\multicolumn{8}{|c|}{Gains (\%)} \\
\hline
$Q_{b,2s}^{(j)}$ vs. $K_{bb,l}^{\text{neigh}(j)}$ & -3 & 13 & 10 & 10 & 12 & 15 & 15 \\
\hline
\hline
$Q_{b,2s}^{\text{neigh}(j)}$ vs. NKS & 30 & 25 & 20 & 18 & 19 & 24 & 16 \\
\hline
\end{tabular}
\end{center}
\caption{Bi-material beam, coupling with SRKS: Krylov cumulated iterations over load increments}
\label{tab:iter_bimat_SRKS}
\end{table}

\begin{table}[!ht]
\begin{center}
\begin{tabular}{|c|cccc:c|}
\hline
\multicolumn{6}{|c|}{Krylov} \\
\hline
\multicolumn{5}{|c:}{with SRKS} & wo SRKS \\
\hline
load inc. & 0.6 & 1 & 1.3 & 1.5 & 1.5 \\
\hline
$S_t^{(\overline{j})}$ & 48 & 164 & 304 & 445 & 481 \\
$K_{bb,l}^{\text{neigh}(j)}$ & 61 & 201 & 351 & 506 & 548 \\
$Q_{b,2s}^{(j)}$ & 47 & 159 & 280 & 410 & 438 \\
\hline
\hline
NKS & 72 & 212 & 362 & 517 & 561 \\
\hline
\multicolumn{6}{c}{ } \\
\hline
\multicolumn{6}{|c|}{Gains (\%) } \\
\hline
$Q_{b,2s}^{(j)}$ vs. $K_{bb,l}^{\text{neigh}(j)}$ & 20 & 20 & 20 & 20 & 20 \\
$Q_{b,2s}^{(j)}$ vs. $S_t^{(\overline{j})}$ & 2 & 3 & 8 & 8 & 9 \\
\hline
\hline
$Q_{b,2s}^{(j)}$ vs. NKS & 35 & 25 & 23 & 21 & 22 \\
\hline
\end{tabular}
\end{center}
\caption{Multiperforated beam, coupling with SRKS: Krylov cumulated iterations over load increments}
\label{tab:iter_multiperf_SRKS}
\end{table}

As expected, SRKS leads to a global decrease of the number of Krylov iterations, observable by comparing the columns ''with'' and ''without'' SRKS of results tables. For the bi-material beam, Krylov iterations are reduced on average by 67\% at last load increment; for multiperforated beam the average reduction is only close to 8\% (a small number of Krylov iterations implies a small number of post-processed Ritz vectors: this could partly explain the less impressive efficiency of SRKS method on this test case). 

Concerning global Newton solver, the cumulated numbers of iterations remained constant with and without SRKS -- as expected, -- they were thus not presented again in this section. 

Tables \ref{tab:iter_bimat_SRKS} and \ref{tab:iter_multiperf_SRKS} confirm observations of previous section. Even if the cumulated numbers of Krylov iterations are decreased thanks to SRKS, the overall gains generated by the new expression $Q_{b,2s}^{(j)}$ remain rather constant, and are even better for bi-material beam (indeed, the classic method NKS suffered from a slight degradation of its overall performance, and the gain of impedance $Q_{b,2s}^{(j)}$ compared to NKS reaches then 24\% at the end of the resolution, in terms of Krylov cumulated iterations). 

\section{Conclusion}

A new approximation of the interface impedance has been developed, in the context of nonlinear substructuring and condensation methods with mixed approach. The expression of the interface impedance introduced here couples both short and long range interactions terms.

The procedure for building such a parameter consists in evaluating, for a given subdomain, the Schur tangent operator of the remainder of the structure (i.e. the optimal value in a linear context), which was originally the best analytic expression we could produce to approximate the optimal interface impedance in the nonlinear context. This evaluation involves a short scale term, basically consisting in the stiffness of the considered subdomain neighbors, and a long scale low rank term, composed of the projection of the Schur tangent operator into the space generated by rigid body modes, thereby capturing long range interactions with the faraway structure. 

Performance of a FETI-2LM solver was studied on a linear case, where the Schur tangent operator of the remainder is exactly the optimal value for the interface impedance -- despite its intractability in practice in parallel resolution processes. Although, as expected, the new additive expression of the impedance did not produce as good results as this optimal value, it managed quite impressive gains, in particular compared to the classical choice made in this framework -- i.e. the stiffness assembled over the neighbors of a subdomain. 

Performance of the mixed nonlinear substructuring and condensation method was also studied, on a plasticity case. Not only the exact computation of Schur tangent operator is not affordable in the framework of parallel distributed computations -- unlike the new expression we build, which was chosen to be inexpensively calculable in parallel, -- but it also was found to achieve not as good results as this new expression. This suggests that the level of accuracy obtained on the representation of a substructure environment with the additive expression of the interface impedance introduced here is increased. 

Eventually, a study of the coupling of the resolution process with a selective reuse procedure of Krylov solver Ritz vectors (SRKS) tends to assess that performance of this new expression is maintained while numbers of Krylov iterations are decreased.
All these considerations are rather promising for implementations at larger scales.

\bibliographystyle{ieeetr}
\bibliography{refs_biblio}

\begin{thebibliography}{10}

\bibitem{crisfield1979faster}
M.~Crisfield, ``{A} faster modified {N}ewton-{R}aphson iteration,'' {\em
  Computer Methods in Applied Mechanics and Engineering}, vol.~20, no.~3,
  pp.~267--278, 1979.

\bibitem{deuflhard1975relaxation}
P.~Deuflhard, ``{A} relaxation strategy for the modified {N}ewton method,'' in
  {\em Optimization and optimal control}, pp.~59--73, Springer, 1975.

\bibitem{dennis1977quasi}
J.~E. Dennis, Jr and J.~J. Mor{\'e}, ``{Q}uasi-{N}ewton methods, motivation and
  theory,'' {\em SIAM review}, vol.~19, no.~1, pp.~46--89, 1977.

\bibitem{zhang1982modified}
L.~Zhang and D.~Owen, ``{A} modified secant {N}ewton method for non-linear
  problems,'' {\em Computers \& Structures}, vol.~15, no.~5, pp.~543--547,
  1982.

\bibitem{mandel1993balancing}
J.~Mandel, ``{B}alancing domain decomposition,'' {\em Communications in
  numerical methods in engineering}, vol.~9, no.~3, pp.~233--241, 1993.

\bibitem{le1994domain}
P.~Le~Tallec, ``{D}omain decomposition methods in computational mechanics,''
  {\em Computational mechanics advances}, vol.~1, no.~2, pp.~121--220, 1994.

\bibitem{rixen1999simple}
D.~J. Rixen and C.~Farhat, ``{A} simple and efficient extension of a class of
  substructure based preconditioners to heterogeneous structural mechanics
  problems,'' {\em International Journal for Numerical Methods in Engineering},
  vol.~44, no.~4, pp.~489--516, 1999.

\bibitem{farhat2001feti}
C.~Farhat, M.~Lesoinne, P.~LeTallec, K.~Pierson, and D.~Rixen, ``{FETI-DP}: a
  dual--primal unified {FETI} method—part {I}: {A} faster alternative to the
  two-level {FETI} method,'' {\em International journal for numerical methods
  in engineering}, vol.~50, no.~7, pp.~1523--1544, 2001.

\bibitem{gosselet2006non}
P.~Gosselet and C.~Rey, ``{N}on-overlapping domain decomposition methods in
  structural mechanics,'' {\em Archives of computational methods in
  engineering}, vol.~13, no.~4, pp.~515--572, 2006.

\bibitem{cresta2007nonlinear}
P.~Cresta, O.~Allix, C.~Rey, and S.~Guinard, ``{N}onlinear localization
  strategies for domain decomposition methods: {A}pplication to post-buckling
  analyses,'' {\em Computer Methods in Applied Mechanics and Engineering},
  vol.~196, no.~8, pp.~1436--1446, 2007.

\bibitem{pebrel2008nonlinear}
J.~Pebrel, C.~Rey, and P.~Gosselet, ``A nonlinear dual domain decomposition
  method: application to structural problems with damage,'' {\em International
  Journal for Multiscale Computational Engineering}, vol.~6, no.~3,
  pp.~251--262, 2008.

\bibitem{bordeu2009balancing}
F.~Bordeu, P.-A. Boucard, and P.~Gosselet, ``{B}alancing domain decomposition
  with nonlinear relocalization: {P}arallel implementation for laminates,'' in
  {\em First international conference on parallel, distributed and grid
  computing for engineering}, pp.~CCP--90, 2009.

\bibitem{negrello2016substructured}
C.~Negrello, P.~Gosselet, C.~Rey, and J.~Pebrel, ``{S}ubstructured formulations
  of nonlinear structure problems--influence of the interface condition,'' {\em
  International Journal for Numerical Methods in Engineering}, 2016.

\bibitem{cresta2008decomposition}
P.~Cresta, {\em {D}{\'e}composition de domaine et strat{\'e}gies de
  relocalisation non-lin{\'e}aire pour la simulation de grandes structures
  raidies avec flambage local}.
\newblock PhD thesis, {\'E}cole normale sup{\'e}rieure de Cachan-ENS Cachan,
  2008.

\bibitem{pebrel2009etude}
J.~Pebrel, P.~Gosselet, and C.~Rey, ``{E}tude du choix des conditions
  d'interface pour des strat{\'e}gies non lin{\'e}aire de d{\'e}composition de
  domaine,'' in {\em Neuvi{\`e}me colloque national en calcul des structures},
  vol.~2, pp.~393--398, 2009.

\bibitem{hinojosa2014domain}
J.~Hinojosa, O.~Allix, P.-A. Guidault, and P.~Cresta, ``{D}omain decomposition
  methods with nonlinear localization for the buckling and post-buckling
  analyses of large structures,'' {\em Advances in Engineering Software},
  vol.~70, pp.~13--24, 2014.

\bibitem{lions1990schwarz}
P.-L. Lions, ``{O}n the {S}chwarz alternating method. {III}: {A} variant for
  nonoverlapping subdomains,'' in {\em Third international symposium on domain
  decomposition methods for partial differential equations}, vol.~6,
  pp.~202--223, SIAM Philadelphia, PA, 1990.

\bibitem{roux2009feti}
F.-X. Roux, ``{A} {FETI-2LM} method for non-matching grids,'' in {\em Domain
  Decomposition Methods in Science and Engineering XVIII}, pp.~121--128,
  Springer, 2009.

\bibitem{klawonn2017new}
A.~Klawonn, M.~Lanser, O.~Rheinbach, and M.~Uran, ``New nonlinear feti-dp
  methods based on a partial nonlinear elimination of variables,'' in {\em
  Domain Decomposition Methods in Science and Engineering XXIII}, pp.~207--215,
  Springer, 2017.

\bibitem{magoules2004optimal}
F.~Magoul{\`e}s, F.-X. Roux, and S.~Salmon, ``Optimal discrete transmission
  conditions for a nonoverlapping domain decomposition method for the helmholtz
  equation,'' {\em SIAM Journal on Scientific Computing}, vol.~25, no.~5,
  pp.~1497--1515, 2004.

\bibitem{gander2011optimal}
M.~J. Gander and F.~Kwok, ``Optimal interface conditions for an arbitrary
  decomposition into subdomains,'' {\em Domain Decomposition Methods in Science
  and Engineering XIX}, vol.~78, pp.~101--108, 2011.

\bibitem{nier1998remarques}
F.~Nier, ``Remarques sur les algorithmes de d{\'e}composition de domaines,''
  {\em S{\'e}minaire {\'E}quations}, vol.~2054, pp.~1--24, 1998.

\bibitem{gander2006osm}
M.~J. Gander, ``{O}ptimized {S}chwarz {M}ethods,'' {\em SIAM Journal on
  Numerical Analysis}, vol.~44, no.~2, pp.~699--731, 2006.

\bibitem{ladeveze2000micro}
P.~Ladev{\`e}ze and D.~Dureisseix, ``{A} micro/macro approach for parallel
  computing of heterogeneous structures,'' {\em International Journal for
  Computational Civil and Structural Engineering}, vol.~1, pp.~18--28, 2000.

\bibitem{SAAVEDRA.2012.1}
K.~Saavedra, O.~Allix, and P.~Gosselet, ``{O}n a multiscale strategy and its
  optimization for the simulation of combined delamination and buckling,'' {\em
  International Journal for Numerical Methods in Engineering}, vol.~91, no.~7,
  p.~772–798, 2012.

\bibitem{SAAVEDRA.2016.hal.1}
K.~Saavedra, O.~Allix, P.~Gosselet, J.~Hinojosa, and A.~Viard, ``{A}n enhanced
  non-linear multi-scale strategy for the simulation of buckling and
  delamination on 3{D} composite plates,'' {\em submitted to Computer Methods
  in Applied Mechanics and Engineering}, 2016.

\bibitem{Hoang2014}
T.~P. Hoang, C.~Japhet, M.~Kern, and J.~Roberts, ``{V}entcell conditions with
  mixed formulations for flow in porous media,'' in {\em Domain Decomposition
  Methods in Science and Engineering XXII} (T.~Dickopf, M.~Gander, L.~Halpern,
  R.~Krause, and L.~Pavarino, eds.), (Lugano (Switzerland)), pp.~531--540,
  2014.

\bibitem{DESMEURE.2011.3.1}
G.~Desmeure, P.~Gosselet, C.~Rey, and P.~Cresta, ``{E}tude de différentes
  représentations des interefforts dans une stratégie de décomposition de
  domaines mixte,'' in {\em Actes du $10^{eme}$ colloque national en calcul des
  structures}, (Giens (Var, France)), 2011.

\bibitem{magoules_algebraic_2006}
F.~Magoulès, F.~X. Roux, and L.~Series, ``{A}lgebraic approximation of
  {D}irichlet-to-{N}eumann maps for the equations of linear elasticity,'' {\em
  Computer Methods in Applied Mechanics and Engineering}, vol.~195, no.~29-32,
  p.~3742–3759, 2006.

\bibitem{Oumaziz2017}
P.~Oumaziz, P.~Gosselet, P.-A. Boucard, and S.~Guinard, ``A non-invasive
  implementation of a mixed domain decomposition method for frictional contact
  problems,'' {\em Computational Mechanics}, Jul 2017.

\bibitem{cai1999restricted}
X.-C. Cai and M.~Sarkis, ``{A} restricted additive {S}chwarz preconditioner for
  general sparse linear systems,'' {\em Siam journal on scientific computing},
  vol.~21, no.~2, pp.~792--797, 1999.

\bibitem{Far94bis}
C.~Farhat and F.~X. Roux, ``Implicit parallel processing in structural
  mechanics,'' {\em Computational Mechanics Advances}, vol.~2, no.~1,
  pp.~1--124, 1994.
\newblock North-Holland.

\bibitem{nouy03b}
P.~Ladev{\`e}ze and A.~Nouy, ``On a multiscale computational strategy with time
  and space homogenization for structural mechanics,'' {\em Computer Methods in
  Applied Mechanics and Engineering}, vol.~192, pp.~3061--3087, 2003.

\bibitem{ibrahimbegovic2003strong}
A.~Ibrahimbegovi{\'c} and D.~Markovi{\v{c}}, ``{S}trong coupling methods in
  multi-phase and multi-scale modeling of inelastic behavior of heterogeneous
  structures,'' {\em Computer Methods in Applied Mechanics and Engineering},
  vol.~192, no.~28, pp.~3089--3107, 2003.

\bibitem{feyel2000fe}
F.~Feyel and J.-L. Chaboche, ``{FE} 2 multiscale approach for modelling the
  elastoviscoplastic behaviour of long fibre {S}i{C}/{T}i composite
  materials,'' {\em Computer methods in applied mechanics and engineering},
  vol.~183, no.~3, pp.~309--330, 2000.

\bibitem{guidault2007two}
P.~Guidault, O.~Allix, L.~Champaney, and J.~Navarro, ``{A} two-scale approach
  with homogenization for the computation of cracked structures,'' {\em
  Computers \& structures}, vol.~85, no.~17, pp.~1360--1371, 2007.

\bibitem{guidault2008multiscale}
P.-A. Guidault, O.~Allix, L.~Champaney, and C.~Cornuault, ``{A} multiscale
  extended finite element method for crack propagation,'' {\em Computer Methods
  in Applied Mechanics and Engineering}, vol.~197, no.~5, pp.~381--399, 2008.

\bibitem{gendre2011two}
L.~Gendre, O.~Allix, and P.~Gosselet, ``{A} two-scale approximation of the
  {S}chur complement and its use for non-intrusive coupling,'' {\em
  International Journal for Numerical Methods in Engineering}, vol.~87, no.~9,
  pp.~889--905, 2011.

\bibitem{bebendorf2003existence}
M.~Bebendorf and W.~Hackbusch, ``Existence of {H}-matrix approximants to the
  inverse {FE}-matrix of elliptic operators with {L}8-coefficients,'' {\em
  Numerische Mathematik}, vol.~95, no.~1, pp.~1--28, 2003.

\bibitem{amestoy2016complexity}
P.~Amestoy, A.~Buttari, J.-Y. L'Excellent, and T.~Mary, ``On the complexity of
  the block low-rank multifrontal factorization,'' {\em Methods and Algorithms
  for Scientific Computing}, vol.~39, pp.~A1710--A1740, 2017.

\bibitem{oumaziz2}
P.~Oumaziz, P.~Gosselet, P.-A. Boucard, and M.~Abbas, ``A parallel non-invasive
  multi-scale strategy for a mixed domain decomposition method with frictional
  contact,'' {\em in preparation}, 2017.

\bibitem{klawonn2001feti}
A.~Klawonn and O.~B. Widlund, ``{FETI} and {N}eumann-{N}eumann iterative
  substructuring methods: connections and new results,'' {\em Communications on
  pure and applied Mathematics}, vol.~54, no.~1, pp.~57--90, 2001.

\bibitem{SPILLANE:2013:FETI_GenEO_IJNME}
N.~Spillane and D.~J. Rixen, ``{Automatic spectral coarse spaces for robust
  {FETI} and {BDD} algorithms},'' {\em Internat. J. Num. Meth. Engin.},
  vol.~95, no.~11, pp.~953--990, 2013.

\bibitem{gosselet2013total}
P.~Gosselet, C.~Rey, and J.~Pebrel, ``{T}otal and selective reuse of {K}rylov
  subspaces for the resolution of sequences of nonlinear structural problems,''
  {\em International Journal for Numerical Methods in Engineering}, vol.~94,
  no.~1, pp.~60--83, 2013.

\bibitem{dubois2012optimized}
O.~Dubois, M.~J. Gander, S.~Loisel, A.~St-Cyr, and D.~B. Szyld, ``The optimized
  schwarz method with a coarse grid correction,'' {\em SIAM Journal on
  Scientific Computing}, vol.~34, no.~1, pp.~A421--A458, 2012.

\bibitem{haferssas2015robust}
R.~Haferssas, P.~Jolivet, and F.~Nataf, ``A robust coarse space for optimized
  schwarz methods: Soras-geneo-2,'' {\em Comptes Rendus Mathematique},
  vol.~353, no.~10, pp.~959--963, 2015.

\bibitem{loisel2015optimized}
S.~Loisel, H.~Nguyen, and R.~Scheichl, ``Optimized schwarz and 2-lagrange
  multiplier methods for multiscale elliptic pdes,'' {\em SIAM Journal on
  Scientific Computing}, vol.~37, no.~6, pp.~A2896--A2923, 2015.

\end{thebibliography}

\end{document}